\documentclass{article}
\usepackage{latexsym}
\usepackage{amsmath,amssymb}
\usepackage[dvips]{graphicx, color}

\newtheorem{thm}{Theorem}

\newtheorem{prop}[thm]{Proposition}

\newtheorem{lem}[thm]{Lemma}
\newtheorem{conj}[thm]{Conjecture} 
\newtheorem{slem}[thm]{Sublemma} 

\newcommand{\Aut}{\mathrm{Aut}}

\newcommand{\het}{\mathrm{ht}}

\newcommand{\z}{\widehat\zeta_X}

\newcommand{\W}{\frak{W}}

\begin{document}
\title{\bf\Large{\ \ Special Uniformity of Zeta Functions}
\newline { I. Geometric Aspect}}
\author{\bf{\large{Lin WENG}}}
\date{}
\maketitle

\noindent
{\it -- In memory of those who lost their lives in the 2011 Tohoku-Kanto Earthquake}
\vskip 0.30cm
\noindent
{\footnotesize{{\bf Abstract:} The special uniformity of zeta functions claims that pure non-abelian
zeta functions coincide with group zeta functions associated to the special linear groups. Naturally associated are three aspects, namely, the analytic, arithmetic, and  geometric aspects. In the first paper of this series, we expose intrinsic geometric structures of our zetas by counting semi-stable  bundles on curves defined over finite fields in terms of their automorphism groups and global sections. We show that such a counting maybe read from Artin zetas which are abelian in nature.
This paper also contains an appendix written by H. Yoshida, one of the driving forces for us to seek group zetas. In this appendix, Yoshida introduces a new zeta as a function field analogue of the group zeta for $SL_2$ for number fields and establishes the Riemann Hypothesis for it.}}

\section{Introduction}
Let $X$ be an irreducible, reduced and regular projective curve of genus $g$ defined over finite field $\mathbb F_q$. Denote its degree $d$  Picard variety by $\mathrm{Pic}^d(X)$. 
This Picard variety admits a natural Brill-Noether stratification defined using the global sections of its elements. Namely,
$$W_X^{\geq i}(d):=\Big\{L\in \mathrm{Pic}^d(X):h^0(X,L)\geq i\Big\}.$$
Set $$W_X^i(d):=W_X^{\geq i}(d)\,\backslash\, W_X^{\geq (i+1)}(d)\qquad\mathrm{and}\qquad
w_X(d;i):=|W_X^i(d)|.$$ 
By the duality and the Clifford lemma, we have for $0\leq d\leq 2g-2$
$$w_X(d,i)=w_X(2g-2-d,i-d+(g-1)),\qquad\mathrm{and}\quad w_X(d,j)=0\ \mathrm{if}\ j\geq \frac{d}{2}+1.$$ 
More generally, the complexity of  $\mathrm{Pic}^d(X)$, or better of $X$, may be measured by $\gamma$ invariants defined as
$$\gamma_X(d):=\frac{1}{q-1}\sum_{i\geq 0}q^i\cdot w_X(d,i).$$
For example, $\gamma_X(0)=\frac{q}{q-1}$, since
$W_X(0)^{\geq 1}=\{\mathcal O_X\};$ and by the vanishing theorem,
$\gamma_X(d)=\frac{N}{q-1}\cdot q^{d-(g-1)},$ if  $d> 2g-2$, where $N=\#X(\mathbb F_q)$.

The aim of this paper is to study the invariants $\gamma$ systematically for semi-stable bundles (so line bundles are included).  In doing so, the first difficulty we face is that: not only the stratification induced by the Brill-Noether loci (so $h^0$ naturally appears) and the automorphism groups (so the factor $q-1$ above is explained) should be taken into the consideration, a refined structure defined by the
 associated poly-stable Jordan-H\"older graded bundles should be treated properly. It is for this purpose  that we start the paper with a construction of
what we call the {\it fat} moduli space ${\bf M}_{X,r}(d)$
of semi-stable bundles. With this moduli space introduced, similarly, then we get
the associated Brill-Noether loci
$W_{X,r}^{\geq i}(d)$ and $W^i_{X,r}(d)$.  Clearly, arithmetic and geometric complexities of $X$ are measured by these $W_{X,r}^{\geq i}(d)$'s.

An effective way to study ${\bf M}_{X,r}(d)$ is to introduce  globally defined  invariants $$\begin{aligned}\alpha_{X,r}(d):=&\sum_{i\geq 0}(q^i-1)\cdot \sum_{\mathcal E\in W_{X,r}^i(d)}\frac{1}{\mathrm{Aut}\,\mathcal E},\\
\beta_{X,r}(d):=&\sum_{i\geq 0} \sum_{\mathcal E\in W_{X,r}^i(d)}\frac{1}{\mathrm{Aut}\,\mathcal E},\\
\gamma_{X,r}(d):=&\sum_{i\geq 0}q^i\cdot \sum_{\mathcal E\in W_{X,r}^i(d)}\frac{1}{\mathrm{Aut}\,\mathcal E}
\end{aligned}$$
 The study of the invariant $\beta$, a classical theme, has a long history. The highest points are the works of Siegel-Weil (see e.g. [HN]), Harder-Narasimhan ([HN]), Desale-Ramanan ([DR]) and Zagier ([Z]). On the other hand, the study of the invariants $\gamma$, or the same $\alpha$, is just about to begin, even it was initiated in our earlier papers on zetas (see e.g., [W1] and the references listed there). The main reason for being less progressed is that, except the fact they fit together very well in our non-abelian zetas, no other uniform structure has been detected.

However, this has been changed dramatically since our work [W3]. In our paper [W3], itself motivated by an old work of Drinfeld on counting rank two 
cuspidal representations associated to curves, we introduce the pure non-abelian zeta functions by counting {\it only those} fat moduli spaces {\it whose degrees are multiples of the ranks}. 
This new idea used in the construction proves to be very crucial: For new pure zetas, not only the associated Riemann Hypothesis can be expected, but a special uniformity  holds. As a direct consequence, we
then have a uniform control of the $\alpha$'s in the case when $d\in r\cdot\mathbb Z$.

To explain this, we need to recall the so-called group zeta functions for function fields ([W3]). Built up on the Artin zetas, these group zeta functions depend only on the group structures, or better, the Lie structures involved. Simply put, the special uniformity claims that our rank $r$ pure zeta function coincides with the
group zeta function associated to the special linear group $SL_r$. A proof of this special uniformity, very much involved, could be obtained using the tools, methods  and constructions
appeared in the studies of trace formulas, such as
the theory of Arthur's analytic truncation and Lafforgue's geo-arithmetic truncation, Langlands-Morris-Siegel's theory of Eisenstain seires, an advanced version of Rankin-Selberg \& Zagier methods, etc. This will be done in [W4], motivated by our earlier work on number fields ([W2]).

In this paper, we will use this special uniformity to investigate the invariants $\alpha$'s and hence also $\gamma$'s. Our final result claims that: For all ranks $r$, the invariants $\alpha_{X,r}(m\cdot r),\ m\in\mathbb Z,$
can be expressed  in terms of the elementary symmetric polynomials of the abelian Weil roots, or equivalently, in terms of 
the coefficients, of the numerator polynomial of the standard Artin zeta function. 
For details, please refer to Prop 2, resp. Prop 6, resp. Thm 7, for line bundles, resp.  higher rank bundles, resp.  rank two bundles.

The study here is expected to play a central role in understanding non-abelian arithmetic of curves.

This paper has an appendix written by H. Yoshida on a new type of zeta functions  and their Riemann Hypothesis. Motivated by our works on zeta functions for number fields associated to $SL_2$ 
and the works of Lagarias-Suzuki  and Ki, Yoshida
introduces a new type zeta functions and shows that these zetas satisfy the Riemann Hypothesis. This work is one of the driving forces for us to search for group zetas for function fields. In fact, we will show that Yoshida's new zetas are essentially the group zetas associated to $SL_2$. Moreover,  the RH established offers an excellent control of the invariants
$\alpha$'s and $\beta$ for rank two bundles using Thm 7.

\section{Special Uniformity of Zetas: Initial State}
\subsection{Refined Geometry of Curves}
Let $X$ be an irreducible, reduced and regular projective curve of genus $g$ defined over a finite field $\mathbb F_q$. For a fixed pair $(r,d)\in\mathbb Z_{>0}\times\mathbb Z$, we then have a naturally associated  {\it fat} moduli space ${\bf M}_{X,r}(d)$ of semi-stable bundles of rank $r$ and degree $d$ defined over $X/\mathbb F_q$.

Indeed, if we denote by $\mathcal M_{X,r}(d)$ the standard moduli space of  semi-stable bundles of rank $r$ and  degree $d$ for $X/\mathbb F_q$. Then
$\mathcal M_{X,r}(d)$ consist of certain equivalence classes $[\mathcal E]$
of semi-stable $\mathbb F_q$-rational bundles $\mathcal E$ of rank $r$ and degree $d$ defined over $X$: Being semi-stable, $\mathcal E$ admits a Jordan-H\"older filtration 
$$0=\mathrm{Fil}^0\mathcal E\subset \mathrm{Fil}^1\mathcal E\subset\dots\subset \mathrm{Fil}^s\mathcal E=\mathcal E$$ of sub-bundles
defined over $\overline X:=X\times_{\mathbb F_q}\overline{\mathbb F_q}$
satisfying $$\mathrm{Gr}^i\mathcal E:=\mathrm{Fil}^i\mathcal E/\mathrm{Fil}^{i-1}\mathcal E,\qquad i=1,2,\dots,r$$ are stable bundles of slope $\frac{d}{r}$. One checks that, while the filtration
is not unique, its associated graded bundle $$\mathrm{Gr}_{\mathrm{JH}}\mathcal E:=\oplus_{i=1}^s\mathrm{Gr}^i\mathcal E$$ is unique and $\mathbb F_q$-rational. By definition, $\mathcal E$ and $\mathcal E'$ are called {\it equivalent} if $\mathrm{Gr}_{\mathrm{JH}}\mathcal E\simeq \mathrm{Gr}_{\mathrm{JH}}\mathcal E'$ as $\mathbb F_q$-rational bundles. 
Denote by $[\mathcal E]$ the equivalence class associated to the $\mathbb F_q$-rational bundle $\mathcal E$. Then we know that the moduli space $\mathcal M_{X,r}(d)$  consists of these equivalence classes of $\mathbb F_q$-rational semi-stable bundles of rank $r$ and  degree $d$ defined over $X$. With these said, by definition, the {\it fat moduli space}  ${\bf M}_{X,r}(d)$ is the space built from  $\mathcal M_{X,r}(d)$ with the point $[\mathcal E]$ replaced by the collection of semi-stable bundles in $[\mathcal E]$, namely, the set $\big\{\mathcal E:\mathcal E\in [\mathcal E]\big\}$ is added at the point $[\mathcal E]$. 

A natural question is to count these 
$\mathbb F_q$-rational semi-stable bundles $\mathcal E$. For this purpose, two invariants, namely, the automorphism group $\Aut(X,\mathcal E)$ and its global sections $h^0(X,\mathcal E)$ can be naturally used.
This then leads to the refined Brill-Noether loci $$W_{X,r}^{\geq i}(d):=\Big\{\mathcal E\in {\bf M}_{X,r}(d):
\min_{\mathcal E\in[\mathcal E]}:h^0(X,\mathcal E)\geq i\Big\}$$ and $$[\mathcal E]^j:=\{\mathcal E\in[\mathcal E]: \dim_{\mathbb F_q}\Aut\,{\mathcal E}\geq j\}.$$ 

To understand the structures of $W_{X,r}^{\geq i}(d)$, as a staring point, we use the following general principles:

\noindent
(1)  The duality: There exist natural
isomorphisms $${\bf M}_{X,r}(d)\to {\bf M}_{X,r}(d+rm),\qquad
\mathcal E\mapsto A^m\otimes\mathcal E$$ and $${\bf M}_{X,r}(d)\to {\bf M}_{X,r}(-d+r(2g-2)),\qquad
\mathcal E\mapsto K_X\otimes\mathcal E^\vee,$$ where $A$ is an Artin line bundle of degree one and $K_X$ denotes the dualizing bundle of $X/\mathbb F_q$.  In particular, $$h^1(X,\mathcal E)=h^0(X,K_X\otimes\mathcal E^\vee).$$

\noindent
(2) The  Riemann-Roch theorem:
$$h^0(X,\mathcal E)-h^1(X,\mathcal E)=\mathrm{deg}(\mathcal E)-\mathrm{rk}(\mathcal E)\cdot (g-1).$$

\noindent
(3)  The vanishing theorem: For semi-stable bundles $\mathcal E$,
$h^1(X,\mathcal E)=0$ if $d(\mathrm E)\geq r(2g-2)+1$. And

\noindent
(4) The Clifford lemma: For semi-stable bundles $\mathcal E$,
$h^0(X,\mathcal E)\leq \mathrm{rk}(\mathcal E)+\frac{\mathrm{deg}(\mathcal E)}{2}$ if $0\leq \mu(\mathcal E)\leq 2g-2$. 

Thus, we only need to understand the moduli spaces ${\bf M}_{X,r}(d)$ with $0\leq d\leq r(g-1),$ 
or better, the refined Brill-Noether loci
$$\Big\{W_{X,r}^{\geq i}(d)\ :\ 0\leq d\leq r(g-1),\ \ 0\leq i\leq r+\frac{d}{2}\Big\}.$$

To effectively investigate these spaces globally, we introduce invariants $\alpha$, $\beta$ and $\gamma$:
$$\begin{aligned}\alpha_{X,r}(d):=&\sum_{\mathcal E\in {\bf M}_{X,r}(d)}\frac{q^{h^0(X,\mathcal E)}-1}{\#\Aut(\mathcal E)},\\
 \beta_{X,r}(d):=&\sum_{\mathcal E\in {\bf M}_{X,r}(d)}\frac{1}{\#\Aut(\mathcal E)},\\
\gamma_{X,r}(d):=&\sum_{\mathcal E\in {\bf M}_{X,r}(d)}\frac{q^{h^0(X,\mathcal E)}}{\#\Aut(\mathcal E)}.\end{aligned}$$  Clearly, $\gamma=\alpha+\beta.$
As $\beta$, a  classical invariant ([HN]), is well-known, so it suffices to study $\alpha$.

\subsection{Pure Non-Abelian Zetas}
The invariants $\alpha$ and $\beta$, hence the geometry of the curve $X$, are systematically
dominated by pure non-abelian zeta functions introduced in [W3]. Recall that these pure zetas, a natural generalization of the Artin zeta function, are
defined as follows:
$$\begin{aligned}\zeta_{X,r}(s):=&\sum_{m=0}^\infty
\sum_{V\in {\bf  M}_{X,r}(d), d=rm}\frac{q^{h^0(X,V)}-1}{\#\mathrm{Aut}(V)}\cdot (q^{-s})^{d(V)},\\\widehat\zeta_{X,r}(s):=&\sum_{m=0}^\infty
\sum_{V\in {\bf  M}_{X,r}(d), d=rm}\frac{q^{h^0(X,V)}-1}{\#\mathrm{Aut}(V)}\cdot (q^{-s})^{\chi(X,V)}.\end{aligned}$$ As usual, set
$$\qquad Z_{X,r}(t):=\zeta_{X,r}(s)\qquad\mathrm{and}\qquad \widehat Z_{X,r}(t):=\widehat\zeta_{X,r}(s)\qquad\mathrm{with}\qquad t:=q^{-s}.$$

\begin{thm} ([W3]) ({\bf Zeta Facts}) (i) $\zeta_{X,1}(s)=\zeta_X(s)$, the Artin zeta function for $X/\mathbb F_q$;

\noindent
(ii) ({\bf Rationality}) There exists a degree $2g$ polynomial $P_{X,r}(T)\in\mathbb Q[T]$ of $T$
such that
$$Z_{X,r}(t)=\frac{P_{X,r}(T)}{(1-T)(1-QT)}\quad\mathrm{with}\quad T=t^r,\ Q=q^r;$$
\noindent
(iii) ({\bf Functional equation}) $$\widehat Z_{X,r}(\frac{1}{qt})=\widehat Z_{X,r}(t),$$

\end{thm}
Indeed, from [W3], we know that
$$\begin{aligned}Z_{X,r}&(t)
=\sum_{m=0}^{(g-1)-1}\alpha_{X,r}(mr)
\cdot \Big(T^{m}+ Q^{(g-1)-m}\cdot T^{2(g-1)-m}\Big)\\
&+\alpha_{X,r}\big(r(g-1)\big)\cdot T^{g-1}+(Q-1)\beta_{X,r}(0)\cdot \frac{T^{g}}{(1-T)(1-QT)}.
\end{aligned}\eqno(1)$$
\subsection{Special Uniformity in Rank One}
In this subsection, we use the uniformity in rank one, namely, Thm 1(i), to determine the invariants $\alpha$ and $\beta$.

Write the Artin zeta for $X$ in Weil's form:
$$\begin{aligned}Z_X(t):=&\exp\Big(\sum_{m=1}^\infty N_m\frac{t^m}{m}\Big)
=:\frac{\sum_{i=0}^{2g}A_it^i}{(1-t)(1- q t)}\\
=&\frac{\prod_{i=1}^{2g}(1-\omega_it)}{(1-t)(1-qt)}
=\frac{\prod_{i=1}^g(1-a_i t+q t^2)}{(1-t)(1- q t)}\\
\end{aligned}$$ where the reciprocal roots $\omega_i, i=1,2,\cdots, 2g$ are arranged in pairs such that $q=\omega_i\cdot\omega_{2g-i}$ and accordingly, $a_i:=\omega_i+\omega_{2g-i}$, $A_0=1$ and $A_i$ are well-known elementary symmetric functions in $\omega_i$'s, and 
$$N_m=q^m+1-\sum_{i=1}^{2g}\omega_i^m=\#X(\mathbb F_{q^m}).$$
Then the uniformity in rank one implies that 
$$\begin{aligned}\alpha_0\cdot\frac{\sum_{i=0}^{2g}A_it^i}{(1-t)(1- q t)}
=&\sum_{m=0}^{(g-1)-1}\alpha_m
\cdot \Big(t^{m}+ q^{(g-1)-m}\cdot t^{2(g-1)-m}\Big)\\
&+\alpha_{g-1}\cdot t^{g-1}+(q-1)\beta_{X}(0)\cdot \frac{t^{g}}{(1-t)(1-qt)}.\end{aligned}$$ Here $$\alpha_m:=\alpha_{X}(m):=\alpha_{X,1}(m), \ \ m=0,1,\cdots,g-1,\ \ \mathrm{and}\qquad\beta_0:=\beta_{X}(0):=\beta_{X,1}(0).$$
This implies that
This implies that $\alpha$'s and $\beta$ satisfy the following system of linear equations
$$\begin{cases}\alpha_0&=\alpha_0\cdot A_0\\
\alpha_1 -(q+1)\alpha_0&=\alpha_0\cdot A_1\\
\alpha_2 -(q+1)\alpha_1+q\alpha_0&=\alpha_0\cdot A_2\\
\cdots\cdots\cdots\cdots\cdots\cdots\cdots\cdots&\cdots\cdots\cdots\\
\alpha_{g-2}-(q+1)\alpha_{g-3}+q\alpha_{g-4}&=\alpha_0\cdot A_{g-2}\\
\alpha_{g-1}-(q+1)\alpha_{g-2}+q\alpha_{g-3}&=\alpha_0\cdot A_{g-1}\\
2q\alpha_{g-2}-(q+1)\alpha_{g-1}+(q-1)\beta_0&=\alpha_0\cdot A_g.\end{cases}$$ Indeed,
writing down  $$\begin{aligned}&\sum_{m=0}^{(g-1)-1}\alpha_m
\cdot \Big(t^{m}+ q^{(g-1)-m}\cdot t^{2(g-1)-m}\Big)(1-t)(1-qt)\\
&+\alpha_{g-1}(1-t)(1-qt)\cdot t^{g-1}+(q-1)\beta_{0}\cdot {t^{g}}\end{aligned}$$
in an explicit form, we get
$$\begin{aligned}\alpha_0&\Big(\sum_{i=1}^{2g}A_it^i\Big)\\
=&\Big(1+ q^{g}t^{2g}\Big)\alpha_0\\
&+\Big(\alpha_1 -(q+1)\alpha_0\Big)\Big(t+q^{g-1}t^{2g-1}\Big)\\
&+\Big(\alpha_2 -(q+1)\alpha_1+q\alpha_0\Big)
\Big(t^2+q^{g-2}t^{2g-2}\Big)\\
&+\cdots\\
&+\Big(\alpha_{g-2}-(q+1)\alpha_{g-3}+q\alpha_{g-4}\Big)\Big(t^{g-2}+q^2t^{g+2}\Big)\\
&+\Big(\alpha_{g-1}-(q+1)\alpha_{g-2}+q\alpha_{g-3}\Big)\Big(t^{g-1}+qt^{g+1}\Big)\\
&+\Big(2q\alpha_{g-2}-(q+1)\alpha_{g-1}+(q-1)\beta_0\Big)\cdot t^{g}\\
\end{aligned}$$
So the above system of equations is obtained by comparing with the coefficients of $t^i$'s.

To go further, note that the inverse matrix of
$$\begin{pmatrix}1&0&0&0&\cdots&0&0&0&0\\
-(q+1)&1&0&0&\cdots&0&0&0&0\\
q&-(q +1)&1&0&\cdots &0&0&0&0\\
0&q &-(q +1)&1&\cdots &0&0&0&0\\
\cdots&\cdots&\cdots&\cdots&\cdots&\cdots&\cdots&\cdots&\\
0&0&0&0&\cdots&1&0&0&0\\
0&0&0&0&\cdots&-(q +1)&1&0&0\\
0&0&0&0&\cdots&q &-(q +1)&1&0\\
0&0&0&0&\cdots&0&q &-(q +1)&1\end{pmatrix}$$
is given by
$$\begin{pmatrix}1&0&0&0&\cdots&0&0&0\\
q +1&1&0&0&\cdots&0&0&0\\
q^2+q +1&q +1&1&0&\cdots &0&0&0\\
q^3+q^2+q +1&q^2+q +1&q +1&1&\cdots &0&0&0\\
\cdots&\cdots&\cdots&\cdots&\cdots&\cdots&\cdots&\\
\frac{q^{g-3}-1}{q -1}&\frac{q^{ (g-4)}-1}{q -1}&\frac{q^{ (g-5)}-1}{q -1}&\frac{q^{ (g-6)}-1}{q -1}&\cdots&1&0&0\\
\frac{q^{ (g-2)}-1}{q -1}&\frac{q^{ (g-3)}-1}{q -1}&\frac{q^{ (g-4)}-1}{q -1}&\frac{q^{ (g-5)}-1}{q -1}&\cdots&q +1&1&0\\
\frac{q^{ (g-1)}-1}{q -1}&\frac{q^{ (g-2)}-1}{q -1}&\frac{q^{ (g-3)}-1}{q -1}&\frac{q^{ (g-4)}-1}{q -1}&\cdots&q^2+q +1&q +1&1\end{pmatrix}.$$
Therefore, $$\begin{cases}\alpha_0&= \alpha_0\\
\alpha_1&=\Big((q+1)+A_1\Big)\cdot \alpha_0\\
\alpha_2&=\Big((q^2+q+1)+(q+1)A_1+A_2\Big)\cdot \alpha_0\\
\alpha_3&=\Big((q^3+q^2+q+1)+(q^2+q+1)A_1+(q+1)A_2+A_3\Big)\cdot \alpha_0\\
\cdots&\cdots\\
\alpha_i&=\Big(\frac{q^{i+1}-1}{q-1}+\frac{q^{i}-1}{q-1}A_1+\cdots+(q+1)A_{i-1}+A_i\Big)\cdot \alpha_0\\
\cdots&\cdots\\
\alpha_{g-2}&=\Big(\frac{q^{g-1}-1}{q-1}+\frac{q^{g-2}-1}{q-1}A_1+\cdots+(q+1)A_{g-3}+A_{g-2}\Big)\cdot \alpha_0\\
\alpha_{g-1}&=\Big(\frac{q^{g}-1}{q-1}+\frac{q^{g-1}-1}{q-1}A_1+\cdots+(q+1)A_{g-2}+A_{g-1}\Big)\cdot \alpha_0\\\end{cases}$$
and
$$2q\alpha_{g-2}-(q+1)\alpha_{g-1}+(q-1)\beta_0=\alpha_0\cdot A_g.$$
Now, by definition, $$\alpha_0=\sum_{L\in\mathrm{Pic}^0(X)}\frac{q^{h^0(X,L)}-1}{q-1}
=\frac{q^{h^0(X,\mathcal O_X)}-1}{q-1}=1$$ and
$$\beta_0=\frac{N}{q-1}\qquad\mathrm{with}\qquad N:=|X(\mathbb F_q)|=\sum_{i=0}^{2g}A_i.$$  This then completes the proof of the following

\begin{prop}  Let $X$ be an  irreducible, reduced regular projective curve $X$ of genus $g$ defined over $\mathbb F_q$.

\noindent
(1) In terms of the Weil's coefficients $A_i$'s of the Artin zeta function, the invariants
 $\alpha$'s and $\beta$ are given by 
$$\begin{cases}\alpha_0&= 1\\
\alpha_1&=(q+1)+A_1\\
\alpha_2&=(q^2+q+1)+(q+1)A_1+A_2\\
\alpha_3&=(q^3+q^2+q+1)+(q^2+q+1)A_1+(q+1)A_2+A_3\\
\cdots&\cdots\\
\alpha_i&=\frac{q^{i+1}-1}{q-1}+\frac{q^{i}-1}{q-1}A_1+\cdots+(q+1)A_{i-1}+A_i\\
\cdots&\cdots\\
\alpha_{g-2}&=\frac{q^{g-1}-1}{q-1}+\frac{q^{g-2}-1}{q-1}A_1+\cdots+(q+1)A_{g-3}+A_{g-2}\\
\alpha_{g-1}&=\frac{q^{g}-1}{q-1}+\frac{q^{g-1}-1}{q-1}A_1+\cdots+(q+1)A_{g-2}+A_{g-1}\\
\beta_0&=\beta_m=\frac{1}{q-1}\Big(A_0+A_1+\cdots+A_{2g}\Big).\end{cases}$$

\noindent
(2) In terms of the invariants $\alpha$'s and $\beta$, the Weil's coefficients $A_i$'s of the Artin zeta function are given by
$$\begin{cases}A_0&=\alpha_0= 1\\
A_1&=\alpha_1 -(q+1)\alpha_0\\
A_2&=\alpha_2 -(q+1)\alpha_1+q\alpha_0\\
\cdots\cdots&\cdots\cdots\cdots\cdots\cdots\cdots\cdots\cdots\cdots\\
A_{g-2}&=\alpha_{g-2}-(q+1)\alpha_{g-3}+q\alpha_{g-4}\\
A_{g-1}&=\alpha_{g-1}-(q+1)\alpha_{g-2}+q\alpha_{g-3}\\
A_g:&=2q\alpha_{g-2}-(q+1)\alpha_{g-1}+(q-1)\beta_0.\end{cases}$$
\end{prop}

\noindent
{\bf Remark.}  The ralation $$2q\alpha_{g-2}-(q+1)\alpha_{g-1}+(q-1)\beta_0=\alpha_0\cdot A_g.$$ implies that
$$\begin{aligned}&2q\Big(\frac{q^{g-1}-1}{q-1}+\frac{q^{g-2}-1}{q-1}A_1+\cdots+(q+1)A_{g-3}+A_{g-2}\Big)\\
&-
(q+1)\Big(\frac{q^{g}-1}{q-1}+\frac{q^{g-1}-1}{q-1}A_1+\cdots+(q+1)A_{g-2}+A_{g-1}\Big)\\
&+\Big(A_0+A_1+\cdots+A_{2g}\Big)=A_g\end{aligned}$$
which can be proved using functional equation directly.

For later use,  we set $$\begin{aligned}\widehat\zeta_X(0):=&\frac{1}{q-1}\Big(A_0+A_1+\cdots+A_{2g}\Big),\\
 \widehat\zeta_X(1):=&\frac{1}{q-1}\Big(A_0q^g+A_1q^{g-1}+\cdots+A_{2g}q^{-g}\Big).\end{aligned}$$
 
 \section{Special Uniformity of Zetas in General}
\subsection{Zetas Associated to $SL_r$}
For $G=SL_r$ with $B$ the standard Borel subgroup consisting of upper triangular matrices, let $T$ be the associated torus consisting of diagonal matrices. Then the root system $\Phi$ associated to $T$ can be realized 
as $$\Phi^+=\Phi_r^+=\{e_i-e_j:1\leq i<j\leq r\}$$ with $\{e_i\}_{i=1}^r$ the standard orthogonal basis of the Euclidean space $V=\mathbb R^r$. It is of type $A_{r-1}$, with  simple roots 
$$\Delta:=\{\alpha_i:=e_i-e_{i+1}:i=1,2,\dots,r-1\},$$
the  Weyl vector $$\rho:=\frac{1}{2}\Big((r-1)e_1+(r-3)e_2+\cdots-(r-3)e_{r-1}-(r-1)e_r\Big),$$ and the Weyl group $W$  the permutation group $S_r$ via the action on the subindex of $e_i$'s. Introduce the corresponding fundamental weights $\lambda_j$'s via
$$\langle\lambda_i,\alpha^\vee_j\rangle=\delta_{ij},\qquad\forall\alpha_j\in \Delta.$$
For each $w\in W$, set $\Phi_w:=\Phi^+\cap w^{-1}\Phi^-$.
For $\lambda\in V_{\mathbb C}$,
introduce then the {\it period of an irreducible, reduced regular projective curve $X/\mathbb F_q$ associated to $SL_r$} by
$$\omega^{SL_r}_X(\lambda):=\sum_{w\in W}\frac{1}{\prod_{\alpha\in\Delta}(1-q^{-\langle w\lambda-\rho,\alpha^\vee\rangle})}\prod_{\alpha\in \Phi_w}
\frac{\z(\langle\lambda,\alpha^\vee\rangle)}{\z(\langle\lambda,\alpha^\vee\rangle+1)}.$$

Corresponding to $\alpha_P=\alpha_{r-1}$, let $$P=P_{r-1,1}=\Big\{\begin{pmatrix}A&B\\ 0&D\end{pmatrix}\in SL_r:A\in GL_{r-1},D\in GL_1\Big\}$$
be the standard parabolic subgroup of $SL_n$ attacted to the partition $(r-1)+1=r$. 
Then the associated root system is given by
$$\Phi_P^+=\Phi_{r-1}^+=\{e_i-e_j:1\leq i<j\leq n-1\},$$ with the simple roots
$$\Delta_P=\Delta_{r-1}=\{\alpha_1,\alpha,\dots,\alpha_{r-2}\},$$ the Weyl vector
 $$\rho_P=\frac{1}{2}\Big((r-2)e_1+(r-4)e_2+\cdots-(r-4)e_{r-2}-(r-2)e_{r-1}\Big),$$ and the Weyl group $$W_P=S_{r-1}\hookrightarrow S_r.$$
As such, the corresponding fundamental weight $\lambda_P$, i.e., that normal to $\Phi_P$, is given by
$$\lambda_P=\lambda_{r-1}=\frac{1}{r}\Big(e_1+e_2+\cdots+e_{r-1}-(r-1)e_{r-1}\Big).$$
For later use, set also $$\W_P:=\{w\in W:\Delta_P\subset w^{-1}(\Delta\cup\Phi^-)\}.$$

Write $$\lambda:=\rho+\sum_{j=1}^{r-1}s_j\lambda_j$$ and set $s:=s_{r-1}$. Then 
 we introduce the {\it period of $X/\mathbb F_q$ for $(SL_r,P)$}
as an one variable function defined by
$$\omega_E^{SL_r/P}(s):=\mathrm{Res}_{s_1=0}\mathrm{Res}_{s_2=0}\cdots
\mathrm{Res}_{s_{r-2}=0}\,\omega^{SL_r}_E(\lambda).$$
This period consists many terms, each of which is a product of certain rational factors of $q^{-s}$ and Atrin zetas. Clear up all zeta factors in the denominators of all terms! The resulting function
is then defined to be the {\it zeta function $\z^{SL_r}(s)$ of $X/\mathbb F_q$ associated to $SL_r$}.
Following the Lie structure exposed in [K], we have the following
\begin{thm}([W3])  ({\bf  Functional Equation})
$$\z^{SL_r}(-r-s)=\z^{SL_r}(s).$$
\end{thm}

In fact, much more can be said: If we set $\Phi_P^+=\Phi_{r-1}^+$, then
the structures of $(\Phi^+\backslash\Phi_P^+)\cap w^{-1}\Phi^{\pm}$ and 
$\Phi_P^+\cap w^{-1}\Phi^\pm,\ (w\in W)$ are rather simple. Consequently, one checks, see e.g., 
a \lq joint work' of Kim and myself, or better, [KKS], that
the minimal factor above eliminating all zeta factors appeared in the denominator of each term of 
$\omega_X^{SL_r/P}(s)$ is given by $\prod_{n=2}^{r-1}\z(n)\cdot\z(s+r)$.
Consequently, we have the following
\begin{prop}  The $SL_r$ zeta function $\z^{SL_r}(s)$ is given by
$$\begin{aligned}\z^{SL_r}(s)=&\omega_X^{SL_r/P}(s)\cdot\prod_{n=2}^{r-1}\z(n)\cdot\z(s+r)\\
=&\sum_{n=1}^rR_n(s)\cdot\z(s+n)\end{aligned}
$$ where
$$R_n(s):=\sum_{\substack{w\in \W_P\\
|(\Phi^+\backslash\Phi_P^+)\cap w^{-1}\Phi^+|=n-1}}C_w\prod_{\alpha\in (w^{-1}\Delta\backslash\Phi_P)}\frac{1}{1-q^{-\langle\lambda_P,\alpha^\vee\rangle s-\het\alpha^\vee+1}}$$ and for $w\in\W_P$,
$$C_w:=\prod_{\alpha\in(w^{-1}\Delta)\cap(\Phi_P\backslash\Delta_P)}\frac{\z(2)^{|\Delta_P^+\cap w^{-1}\phi^+|}}{1-q^{1-\het\alpha^\vee}}\prod_{\alpha\in(\Phi_P^+\backslash\Delta_P)\cap w^{-1}\Phi^+}\frac{\z(\het\alpha^\vee+1)}{\z(\het\alpha^\vee)}.$$
\end{prop}

\subsection{Special Uniformity: Analytic Structure}
Using the tools, methods and constructions in the study of trace formula, such as the theory of Arthur's analytic truncations and Lafforgue's  geo-arithmetic truncation,  Langlands-Siegel's theory of Eisenstein series, 
an advanced version of Rankin-Selberg \& Zagier method, etc..., we, in [W4], will show that the rank $r$ pure zeta function and the $SL_r$ zeta function are essentially the same.
To be more precise, we have the following
\begin{thm} ([W4]) ({\bf Special Uniformity})
For an irreducible, reduced regular projective curve $X$ of genus $g$ defined over $\mathbb F_q$,
$$\widehat\zeta_{X,r}(s)=\alpha_{X,r}(0)\cdot \z^{SL_r}(-rs).$$
\end{thm}

As a direct consequence, we see that
$$\z^{SL_r}(-rs)=\frac{\sum_{i=0}^{2g}A_{X,r}(i)T^i}{(1-T)(1-QT)}\cdot\frac{1}{T^{g-1}},$$
for a certain degree $2g$ polynomial $\sum_{i=0}^{2g}A_{X,r}(i)T^i$ of $T$. Moreover,
by definition, the right hand side can be expressed in terms of Artin zeta function up to some rational factors of $q$ depending only on $SL_r$. Consequently, the coefficients
 $A_{X,r}(i)$'s can be read from Artin zeta function $\z(s)$, and so can be expressed in terms of certain combinations of
elementary symmetric polynomials $A_j$'s of the Weil roots $\omega_k,\, k=1,2,\cdots,2g$,
with the help of certain rational functions of $q$ depending only on $SL_r$. This then leads to the following
\begin{prop} (1) In terms of the $SL_r$ zeta coefficients $A_{X,r}(i)$'s, the invariants $\alpha_{X,r}(i)$'s can be effectively
calculated as follows :
{\footnotesize$$\begin{cases}\alpha_{X,r}(0)&= \alpha_{X,r}(0)\\
\alpha_{X,r}(r)&=\Big((q^r+1)+A_{X,r}(1)\Big)\cdot \alpha_{X,r}(0)\\
\alpha_{X,r}(2r)&=\Big((q^{2r}+q^r+1)+(q^r+1)A_{X,r}(1)+A_{X,r}(2)\Big)\cdot \alpha_{X,r}(0)\\
\alpha_{X,r}(3r)&=\Big((q^{3r}+q^{2r}+q^r+1)+(q^{2r}+q^r+1)A_{X,r}(1)+(q^r+1)A_{X,r}(2)+A_{X,r}(3)\Big)\cdot \alpha_{X,r}(0)\\
\cdots&\cdots\\
\alpha_{X,r}(ri)&=\Big(\frac{q^{(i+1)r}-1}{q^r-1}+\frac{q^{ir}-1}{q^r-1}A_{X,r}(1)+\cdots+(q^r+1)A_{X,r}({i-1})+A_{X,r}(i)\Big)\cdot \alpha_{X,r}(0)\\
\cdots&\cdots\\
\alpha_{X,r}({r(g-2)})&=\Big(\frac{q^{(g-1)r}-1}{q^r-1}+\frac{q^{(g-2)r}-1}{q^r-1}A_{X,r}(1)+\cdots
+(q^r+1)A_{X,r}({g-3})+A_{X,r}({g-2})\Big)\cdot \alpha_{X,r}(0)\\
\alpha_{X,r}({r(g-1)})&=\Big(\frac{q^{gr}-1}{q^r-1}+\frac{q^{(g-1)^r}-1}{q^r-1}A_{X,r}(1)+\cdots+(q^r+1)A_{X,r}({g-2})+A_{X,r}({g-1})\Big)\cdot \alpha_{X,r}(0)\\
\beta_{X,r}(0)&=q^{(g-1)\cdot\frac{r^2-r}{2}}\cdot\sum_{\substack{n_1,\dots,n_s>0,\\ n_1+\cdots+n_k=r}}
\frac{(-1)^{k-1}}{\prod_{j=1}^{k-1}(q^{n_j+n_{j+1}}-1)}\prod_{j=1}^k \prod_{i=1}^{n_j}\widehat\zeta_X(i).\end{cases}$$}

\noindent
(2) In terms of the invariants $\alpha$'s and $\beta$, the $SL_r$ zeta coefficients $A_{X,r}(i)$'s can be calculated as follows
{\footnotesize$$\begin{cases}A_{X,r}(0)&=A_{X,r}(0)\\
A_{X,r}(1)&=\Big(\alpha_{X,r}(r) -(q^r+1)\alpha_{X,r}(0)\Big)A_{X,r}(0)\\
A_{X,r}(2)&=\Big(\alpha_{X,r}(2r) -(q^r+1)\alpha_{X,r}(r)+q^r\alpha_{X,r}(0)\Big)A_{X,r}(0)\\
\cdots\cdots&\cdots\cdots\cdots\cdots\cdots\cdots\cdots\cdots\cdots\\
A_{X,r}({g-2})&=\Big(\alpha_{X,r}({(g-2)r})-(q^r+1)\alpha_{X,r}({(g-3)r})+q^r\alpha_{X,r}({(g-4)r})\Big)A_{X,r}(0)\\
A_{X,r}({g-1})&=\Big(\alpha_{X,r}({(g-1)r})-(q^r+1)\alpha_{X,r}({(g-2)r})+q^r\alpha_{X,r}({(g-3)r})\Big)A_{X,r}(0)\\
A_{X,r}(g)&=\Big(2q^r\alpha_{X,r}({(g-2)r})-(q^r+1)\alpha_{X,r}({(g-1)r})+(q^r-1)\beta_{X,r}(0)\Big)A_{X,r}(0).\end{cases}$$}

\noindent
In particular, the rank $r$  pure zeta functions can be read explicitly from the Artin zeta function.
\end{prop} 

\noindent
{\it Proof.} The result on $\beta_{X,r}(0)$ is obtained from the following
\begin{thm} ([HN], qualitative; [DR], quantitative, but not explicit; 
 Zagier [Z])
$$\beta_{X,r}(d)=\sum_{n_1,\dots,n_s>0,\ \sum n_i=r}q^{(g-1)\sum_{i<j}n_in_j}\prod_{i=1}^{s-1}\frac{q^{(n_i+n_{i+1})\{n_1+\cdots+n_i)d/n\}}}{1-q^{n_i+n_{i+1}}}\cdot \prod_{i=1}^s v_{n_i}(q),$$ where $$v_n(q):=\frac{\prod_{i=1}^{2g}(1-\omega_i)}{q-1}q^{(r^2-1)(g-1)}\zeta_X(2)\cdots\zeta_X(n).$$\end{thm}
As for others, a proof can be obtained from a similar discussion as in \S1.3,  with $q$ replaced by $Q$, etc. We leave the details to the reader

\section{Special Uniformity in Level Two}
Thus to explicitly determine the invariants $\alpha_{X,r}(i)$'s from Artin zeta functions, we need to effectively determine the level $r$
coefficients $A_{X,r}(i)$'s appeared in the $SL_r$ zeta function $$\z^{SL_r}(-rs)=\frac{\sum_{i=0}^{2g}A_{X,r}(i)T^i}{(1-T)(1-QT)}\cdot\frac{1}{T^{g-1}}.$$ For general level, this still proves to be  complicated. But for the level two, we know that
$$\begin{aligned}\widehat\zeta_{X,2}(-2s)=&\frac{\widehat\zeta_X(2s)}{1-q^{-2s+2}}+\frac{\widehat\zeta_X(2s-1)}{1-q^{2s}}\\
=&\frac{\sum_{i=0}^{2g}A_iT^{ i}}{(1-T)(1-qT)}\cdot \frac{1}{T^{g-1}}\cdot\frac{1}{1-q^2T}+\frac{\sum_{i=0}^{2g}A_iq^iT^i}{(1-qT)(1-q^2T)}\cdot \frac{1}{(qT)^{g-1}}\cdot\frac{1}{1-\frac{1}{T}}\\
=&
\frac{q^{g-1}\sum_{i=0}^{2g}A_iT^{ i}-
T\sum_{i=0}^{2g}A_iq^iT^i}{(qT)^{g-1}\cdot(1-T)(1-qT)(1-q^2T)}.\end{aligned}.$$
Thus the coefficients $A_{X,2}(i)$'s can be written down precisely. That is, we have the following
\begin{thm} For rank two invariants $\alpha_{X,2}(mr),\ m=0,1,2,\cdots,g-1$,  we have
{\footnotesize$$\begin{cases}\alpha_{X,2}(0)&= q^{g-1}\cdot \widehat\zeta_X(1)\\
\alpha_{X,2}(2)&=\Big((q^2+1)+A_{X,2}(1)\Big)\cdot q^{g-1}\cdot \widehat\zeta_X(1)\\
\alpha_{X,2}(4)&=\Big((q^4+q^2+1)+(q^2+1)A_{X,2}(1)+A_{X,2}(2)\Big)\cdot q^{g-1}\cdot \widehat\zeta_X(1)\\
\alpha_{X,2}(6)&=\Big((q^6+q^4+q^2+1)+(q^4+q^2+1)A_{X,2}(1)+(q^2+1)A_{X,2}(2)+A_{X,2}(3)\Big)\cdot q^{g-1}\cdot \widehat\zeta_X(1)\\
\cdots&\cdots\\
\alpha_{X,2}(2i)&=\Big(\frac{q^{2(i+1)}-1}{q^2-1}+\frac{q^{2(i)}-1}{q^2-1}A_{X,2}(1)+\cdots+(q^2+1)A_{X,2}({i-1})+A_{X,2}(i)\Big)\cdot q^{g-1}\cdot \widehat\zeta_X(1)\\
\cdots&\cdots\\
\alpha_{X,2}({2(g-2)})&=\Big(\frac{q^{2(g-1)}-1}{q^2-1}+\frac{q^{2(g-2)}-1}{q^2-1}A_{X,2}(1)+\cdots+(q^2+1)A_{X,2}({g-3})+A_{X,2}({g-2})\Big)\cdot q^{g-1}\cdot \widehat\zeta_X(1)\\
\alpha_{X,2}({2(g-1)})&=\Big(\frac{q^{2g}-1}{q^2-1}+\frac{q^{2(g-1)}-1}{q^2-1}A_{X,2}(1)+\cdots+(q^2+1)A_{X,2}({g-2})+A_{X,2}({g-1})\Big)\cdot q^{g-1}\cdot \widehat\zeta_X(1)\\
\beta_{X,2}(0)
&=q^{2(g-1)}\cdot\Big(\z(2)-\frac{1}{q^2-1}\cdot\z(1)^2 \Big)\\
\end{cases}$$}
where
{\footnotesize $$\begin{cases}A_{X,2}(1)&=q^{-(g-1)}\Big(a_1q^{g-1}+a_0(q^g-1)\Big)\\
A_{X,2}(2)&=q^{-(g-2)}\Big(a_2q^{g-2}+a_1\Big(q^{g-1}-1\Big)+a_0(q^g-1)\Big)\\
\cdots&\cdots\\
A_{X,2}({g-i})&=q^{-i}\Big(a_{g-i}q^i+a_{g-(i+1)}(q^{i+1}-1)+\cdots+a_2(q^{g-2}-1)+a_1(a^{g-1}-1)+a_0(q^g-1)\Big)\\
\cdots&\cdots\\
A_{X,2}({g-2})&=q^{-2}\Big(a_{g-2}q^2+a_{g-3}(q^3-1)+\cdots+a_2(q^{g-2}-1)+a_1(a^{g-1}-1)+a_0(q^g-1)\Big)\\
A_{X,2}({g-1})&=q^{-1}\Big(a_{g-1}q+a_{g-2}(q^2-1)+\cdots+a_2(q^{g-2}-1)+a_1(a^{g-1}-1)+a_0(q^g-1)\Big)\\
A_{X,2}(g)&=\Big(a_g+a_{g-1}(q-1)+a_{g-2}(q^2-1)+\cdots+a_2(q^{g-2}-1)+a_1(a^{g-1}-1)+a_0(q^g-1)\Big)\end{cases}$$}

\end{thm}

\noindent
{\it Proof.} We begin with the following consequence of the functional equation.

\begin{lem} 
$$\begin{aligned}(1)&\quad q^{g-1}\sum_{i=0}^{2g}A_iT^{i}-
T\sum_{i=0}^{2g}A_iq^iT^i=(1-qT)\cdot\Big[\sum_{i=0}^{g-1}b_i(T^i+T^{2g-i})+b_gT^g\Big].\\
(2)&\quad  \Big(T^{g-1}\cdot(1-T)(1-q^2T)\Big)\cdot\widehat\zeta_{X}^{SL_2}(-2s)\\
=&\sum_{i=0}^{g-1}\frac{A_i}{q^g}\Big[q^{g-i}(qT)^{2g-i}+\Big(q^{g-i}-1\Big)\frac{1-(qT)^{2g-2i-1}}{1-qT}
(qT)^{i+1}+q^{g-i}(qT)^i\Big].\end{aligned}$$
\end{lem}

\noindent
{\it Proof.} By the functional equation for Artin zeta functions, we have
{\footnotesize$$\begin{aligned}&q^{g}\sum_{i=0}^{2g}A_iT^{ i}-
qT\sum_{i=0}^{2g}A_iq^iT^i\\
=&q^{g}\Big[\sum_{i=0}^{g-1}\Big(A_iT^{ i}+A_{2g-i}T^{2g-i}\Big)+A_gT^g\Big]-
qT\Big[\sum_{i=0}^{g-1}\Big(A_iq^iT^i+A_{2g-i}q^{2g-i}T^{2g-i}\Big)+A_gq^gT^g\Big]\\
=&\Big[\sum_{i=0}^{g-1}A_i\Big(q^gT^{ i}+(qT)^{2g-i}\Big)+A_g(qT)^g\Big]-
qT\Big[\sum_{i=0}^{g-1}A_i\Big((qT)^i+q^{g-i}(qT)^{2g-i}\Big)+A_g(qT)^g\Big]\\
=&\sum_{i=0}^{g-1}A_i\Big[\Big(q^{g-i}(qT)^{ i}+(qT)^{2g-i}\Big)-qT\Big((qT)^i+q^{g-i}(qT)^{2g-i}\Big)\Big]+A_g(1-qT)(qT)^g\\
\end{aligned}$$}
 Thus, to prove the lemma, it suffices to prove
 the following elementary
 
\begin{slem}
$$\begin{aligned}&\Big(q^{g-i}x^{i}+x^{2g-i}\Big)-x\Big(x^i+q^{g-i}x^{2g-i}\Big)\\
=&(1-x)\Big[q^{g-i}x^{2g-i}+\Big(q^{g-i}-1\Big)\Big(x^{2g-2i-2}+\cdots+x+1\Big)x^{i+1}+q^{g-i}x^i\Big]\end{aligned}$$
\end{slem}

We leave a proof of this sublemma to the reader.

To continue our proof of the theorem, let us evaluate the coefficients $A_{X,2}(i)$ of $T^i$ in the polynomial
$$\sum_{i=0}^{g-1}\frac{A_i}{q^g}\Big[q^{g-i}(qT)^{2g-i}+\Big(q^{g-i}-1\Big)\frac{1-(qT)^{2g-2i-1}}{1-qT}
(qT)^{i+1}+q^{g-i}(qT)^i\Big].$$
This in practice means  to pin down the coefficients of $X^i$ in the polynomial $$A_g X^g+\sum_{i=0}^{g-1}A_i\Big[q^{g-i}X^{2g-i}+\Big(q^{g-i}-1\Big)\Big(X^{2g-2i-2}+\cdots+X+1\Big)X^{i+1}+q^{g-i}X^i\Big]$$
 where we set $X:=qT$.
 
 By an elementary but dull calculation, we have
{\footnotesize $$\begin{aligned}&a_g X^g+\sum_{i=0}^{g-1}a_i\Big[q^{g-i}X^{2g-i}+\Big(q^{g-i}-1\Big)\Big(X^{2g-2i-2}+\cdots+X+1\Big)X^{i+1}+q^{g-i}X^i\Big]\\
=&a_0q^{g}\Big(X^{2g}+X^{2g-1}+\cdots+1\Big)
\\
&+\Big(a_1q^{g-1}-a_0\Big)\Big(X^{2g-1}+X^{2g-2}+\cdots+X^1\Big)
\\
&+\Big(a_2q^{g-2}-a_1\Big)\Big(X^{2g-2}+X^{2g-3}+\cdots+X^2\Big)
\\
&+\cdots\\
&+\Big(a_{g-2}q^{2}-a_{g-3}\Big)\Big(X^{g+2}+X^{g+1}+X^g+X^{g-1}+X^{g-2}\Big)
\\
&+\Big(a_{g-1}q-a_{g-2}\Big)\Big(X^{g+1}+X^{g}+X^{g-1}\Big)\\
&+\Big(a_g-a_{g-1}\Big)
\Big(X^g\Big)\\
\end{aligned}$$}
Tide all the lose ends up, we then complete the proof of the theorem.

This theorem tells us that the invariants $\alpha_{X,2}(2m)$'s can be effectively
calculated in terms of the Artin zeta function. Thus with the Riemann Hypothesis for Artin zeta functions, we have good controls on $\alpha$'s.
But this is rather remote: after all, the expressions, while explicit, still appear to be very complicated. Thus it is much more crucial to have the following 
\begin{conj}(Riemann Hypothesis)  For all $r\geq 1$,
$$\widehat\zeta_{X,r}(s)=0 \qquad\Rightarrow\qquad\mathrm{Re}(s)=\frac{1}{2}.$$
\end{conj}

In this direction, we have the following

\begin{thm} (Yoshida) ({\bf Riemann Hypothesis})  $$\widehat\zeta_{X}^{SL_2}(s)=0 \qquad\Rightarrow\qquad\mathrm{Re}(s)=\frac{1}{2}.$$
\end{thm}

Indeed, in the appendix, Yoshida introduces the following
new zetas as a functional analogue of  group zeta associated to $SL_2$ for number fields.
$$\widehat\zeta_2(s):=\frac{(1+q^s)\widehat\zeta(2s)}{1-q^{1-s}}
-\frac{q^{-s}(1+q^{1-s})\widehat\zeta_X(2s-1)}{1-q^{-s}}.$$
Moreover he shows the following
\begin{thm}(Yoshida)
$$\widehat\zeta_{2}(s)=0 \qquad\Rightarrow\qquad\mathrm{Re}(s)=\frac{1}{2}.$$
\end{thm}

This discovery was one of the driving forces for us to find the group zetas for function fields.
Now as we know, Yoshida's zeta is essentially the $SL_2$ zeta:
$$\widehat\zeta_2(s)=\Big((1+q^s)(1+q^{1-s})\Big)\cdot\widehat\zeta_X^{SL_2}(s).$$
This then also completes the proof of Thm 12. 
\vskip 0.30cm
We end this paper with the following comments.
Practically, the difficulty of counting semi-stable bundles comes form the fact that direct summands of the associated 
Jordan-H\"older graded bundle, or equivalently, the Jordan-H\"older filtrations, of an 
$\mathbb F_q$-rational semi-stable bundle in general would not be defined over $X/\mathbb F_q$, but
rather its scalar extension $X_n/\mathbb F_{q^n}$. Theoretically, this is the junction point where the abelian and non-abelian ingredients of curves interact. For examples, torsions of Jacobians, Weierstrass points and stable but not absolutely stable bundles are closely related and hence get into the picture naturally. To expose such intrinsic structures is the main theme of our 
study  on arithmetic aspect of the uniformity of zetas. 
\vskip 0.30cm
\centerline{\bf REFERENCES}
\vskip 0.30cm
\noindent
[DR] U.V. Desale \& S. Ramanan, Poincare polynomials of the variety of stable bundles, Math. Ann 26 (1975) 233-244
\vskip 0.10cm
\noindent[HN] G. Harder \& M.S. Narasimhan, On the cohomology groups of moduli spaces of vector bundles on curves, Math. Ann. 212, 215-248 (1975)
\vskip 0.10cm
\noindent 
[KKS] H. Ki, Y. Komori \& M. Suzuki, On the zeros of Weng zeta functions for Chevalley groups, arXiv:1011.4583
\vskip 0.20cm  
\noindent
[K] Y. Komori, Functional equations for Weng's zeta functions for $(G,P)/\mathbb{Q}$, Amer. J. Math., to appear
\vskip 0.10cm
\noindent
[W1] L. Weng, Non-abelian zeta function for function fields, Amer. J. Math., 127 (2005), 973-1017
\vskip 0.10cm
\noindent
[W2] L. Weng, A geometric approach to $L$-functions, in {\it Conference on L-Functions}, pp. 219-370, World Sci (2007)
\vskip 0.10cm
\noindent
[W3] L. Weng, Zeta functions for function fields, preprint, 2012 arXiv:1202.3183
\vskip 0.10cm
\noindent
[W4] L. Weng, Special Uniformity of Zeta Functions, in preparation 
\vskip 0.10cm
\noindent
[Z] D. Zagier, Elementary aspects of the Verlinde formula and the Harder-Narasimhan-Atiyah-Bott formula, in {\it Proceedings of the Hirzebruch 65 Conference on Algebraic Geometry}, 445-462 (1996)
\vskip 0.5cm
\noindent
{\bf Lin WENG}\footnote{
This work is partially supported by JSPS.}

\noindent
Institute for Fundamental Research, The $L$-Academy {\it and}

\noindent
Graduate School of Mathematics, Kyushu University,
 Fukuoka 819-0395

\noindent
E-Mail: weng@math.kyushu-u.ac.jp
\eject

\noindent
{\bf\large{Appendix:}}
\vskip 0.20cm
\qquad{\bf\large{New Zeta Functions and the Riemann Hypothesis}}
\vskip 0.30cm
\hskip 5.0cm{\bf{\large{H. Yoshida}}}
\vskip 1.0cm
\noindent
{\bf \S 1.} 
We start with the following technical lemma.
\vskip 0.20cm
\noindent
{\bf Lemma 1.} {\it Fix  a real number  $q>1$.
Let $\alpha,\,\beta\in\mathbb C$, $\alpha\beta=q$. We put $c=\alpha+\beta$ and assume that $c\in\mathbb R, \ |c|\leq q+1$. Then for $w\in\mathbb C$, we have}
$$|w-\alpha|\,|w-\beta|\,>\,|1-\alpha w|\,|1-\beta w|\qquad \text{if}\  \ |w|<1,\leqno(1)$$
$$|w-\alpha|\,|w-\beta|\,<\,|1-\alpha w|\,|1-\beta w|\qquad \text{if}\  \ |w|>1.\leqno(1')$$

\noindent
{\it Proof.} We have
$$\begin{aligned}&|w-\alpha|^2|w-\beta|^2=(w-\alpha)(w-\beta)(\overline w-\overline\alpha)(\overline w-\overline\beta)\\
=&\Big(w^2-cw+q\Big)\times\Big({\overline w}^2-c\overline w+q\Big)\\
=&|w|^4-c^2|w^2|+q(w^2+{\overline w}^2)-c|w|^2(w+\overline w)+q^2-cq(w+\overline w),\end{aligned}$$
$$\begin{aligned}&|1-\alpha w|^2|1-\beta w|^2=(1-\alpha w)(1-\beta w)(1-\overline \alpha\overline w)(1-\overline \beta\overline w)\\
=&\Big(1-cw+q w^2\Big)\times\Big(1-c{\overline w}+q{\overline w}^2\Big)\\
=&1-c^2|w|^2+q^2|w|^4+q(w^2+{\overline w}^2)-cq|w|^2(w+\overline w)-c(w+\overline w).\end{aligned}$$
Subtract the second formula from the first. Then we obtain
$$\begin{aligned}&(q^2-1)(1-|w|^4)+(q-1)(|w|^2-1)(cw+c\overline w)\\=&(q-1)(1-|w|^2)\Big[(q+1)(1+|w|^2)-(cw+c\overline w)\Big].\end{aligned}$$
Put $r=|w|$. Our inequalities (1) and ($1'$) follow since
$$(q+1)(1+r^2)-2|c|r>0,\qquad r\not=1.$$

Now let $g$ be a positive integer and we consider the function $X_1(s)$ and $X(s)$ defined by

$$X_1(s):=\prod_{i=1}^g(1-\alpha_i q^{-s})(1-\beta_i q^{-s}),\qquad s\in\mathbb C, \leqno(2)$$

$$X(s):=\frac{X_1(s)}{(1-q^{-s})(1-q^{1-s})}, \qquad s\in\mathbb C.\leqno(3)$$ Here $\alpha_i$ and $\beta_i$ are complex numbers. We assume that
$$\alpha_i\beta_i=q,\qquad 1\leq i\leq g.\leqno(4)$$ Then, easily, we have the following
\vskip 0.20cm
\noindent
{\bf Lemma 2.} {\it We have the functional equations}
$$X_1(1-s)=q^{g(2s-1)}X_1(s),\qquad X(1-s)=q^{(g-1)(2s-1)}X(s).$$

We define $$Y(s)=q^{g-1)(s-1/2)}X(s).$$ From Lemma 2, we obtain the functional equation
$$Y(1-s)=Y(s).$$
\vskip 0.30cm
\noindent
{\bf \S2.} Weng's higher rank zeta function for the rational number field, in
the case of rank 2, is equal to
$$\widehat\zeta_{\mathbb Q,2}(s)=\frac{\widehat\zeta(2s)}{s-1}-\frac{\widehat\zeta(2s-1)}{s},$$
where $\widehat\zeta(s):=\pi^{-s/2}\Gamma(s/2)\zeta(s)$ ([W]). The Riemann hypothesis for $\widehat\zeta_{\mathbb Q,2}(s)$ is proved
by Lagarias-Suzuki ([LS]) and independently by Ki ([K]). 

Recall that the zeta function of a projective smooth algebraic
curve of genus $g$ defined over the finite field with $q$ elements has the same
form as $X(s)$. (In the geometric case, we have $\beta_i=\overline{\alpha_i},\ |\alpha_i|=\sqrt q,\ 1\leq i\leq g$.)
We search a function field analogue $\widehat\zeta_2(s)$ of Weng's rank 2 zeta function in
the form $$\widehat\zeta_2(s)=C_1(s)\frac{Y(2s)}{1-q^{1-s}}-C_2(s)\frac{q^{-s}Y(2s-1)}{1-q^{-s}}.$$ Here $C_i(s)$ is a rational function of $q^{-s}$. The functional equation
$$\widehat\zeta_2(1-s)=\widehat\zeta_2(s)$$ holds if and only if
$$C_2(s)=C_1(1-s).\leqno(6)$$
But before we search for the exact form, first, we give a sufficient condition for $\widehat\zeta_2(s)$ to satisfy the Riemann hypothesis.

\vskip 0.20cm
\noindent
{\bf Theorem 3.} {\it Assume that $\alpha_i+\beta_i\in\mathbb R,\ |\alpha_i+\beta_i|\leq q+1$ for $1\leq i\leq g$. We further assume that $C_1(s)$ has the form
$$C_1(s)=q^{as}(1+q^{-s})q^{-hs}\prod_{j=1}^h(1-\gamma_j q^{s-1/2})(1-\delta_j q^{s-1/2}).$$ Here for $1\leq j\leq h$, $\gamma_j$ and $\delta_j$ are complex numbers such that $\gamma_j\delta_j=q,$
$\gamma_j+\delta_j\in\mathbb R$, $|\gamma_j+\delta_j|\leq q+1$ and $a$ is a non-negative real number. Then $\widehat\zeta_2(s)$ satisfies the Riemann hypothesis.}
\vskip 0.20cm
\noindent
{\it Proof.} $Y(s)$ has simple poles at $s$ when $q^s=1,\ q.$
Hence $Y(2s)$ has poles when $q^{2s}=1,\,q$ and $Y(2s-1)$ has poles hen $q^{2s}=q,\,q^2$. We see that $\widehat\zeta_2(s)$
has poles only when $q^{2s} = 1,\, q,\, q^2$. If $q^{2s} = 1$ or $q^2$, 
we see that $s$ is a pole of
$\widehat\zeta_2(s)$. Let $s_0$ be a zero of $\widehat\zeta_2(s)$. Then we have $q^{2s_0}\not=1,\,q^2$.

For simplicity of notation, we write $s_0$ as $s$. As $s$ is a zero of $\widehat\zeta_2(s)$, we
have
$$C_1(s)(1 -q^{-s})Y (2s) = C_2(s)q^{-s}(1 - q^{1-s})Y (2s - 1).\leqno(7)$$
We substitute $Y (2s)$ by $Y (1 - 2s)$ and use the relation
$$Y(s)=q^{(g-1)(s-1/2)}\frac{X_1(s)}{(1-q^{-s})(1-q^{1-s})}.$$
Then (7) is equivalent to
$$\begin{aligned}&C_1(s)q^{(g-1)(1/2-2s)}(1 - q^{-s})\frac{X_1(1-2s)}{(1-q^{2s-1})(1-q^{2s})}\\
=&C_2(s)q^{(g-1)(2s-3/2)}q^{-s}(1 - q^{1-s})\frac{X_1(2s-1)}{(1-q^{1-2s})(1-q^{2-2s})}.\end{aligned}$$
This is equivalent to
$$\begin{aligned}&C_1(s)q^{(g-1)(1/2-2s)}(1 - q^{-s})(1-q^{1-2s})(1-q^{2-2s})
X_1(1 - 2s)\\
=&C_2(s)q^{(g-1)(2s-3/2)}q^{-s}(1 - q^{1-s})(1-q^{2s-1})(1-q^{2s})
X_1(2s-1)
.\end{aligned}\leqno(8)$$
By
$$(1 - q^{2s-1}) = -q^{2s-1}(1 - q^{1-2s}),$$
(8) can be transformed to
$$\begin{aligned}&C_1(s)(1 - q^{-s})(1 - q^{2-2s})X_1(1 - 2s)\\
=&C_2(s)q^{(g-1)(4s-2)}q^{3s-1}(1 - q^{1-s})(1 - q^{-2s})X_1(2s - 1),\end{aligned}$$
which is equivalent to
$$C_1(s)(1 + q^{1-s})X_1(1 - 2s) = C_2(s)q^{(g-1)(4s-2)}q^{3s-1}(1 + q^{-s})X_1(2s - 1).$$

Now put $s = 1/2 + i z$. Then, using (6), (9) can be written as
$$\begin{aligned}&C_1(1/2 + i z)(1 + q^{1/2-iz})X_1(-2 i z)\\
=&C_1(1/2 - iz)q^{(g-1)4iz}q^{1/2+3iz}(1 + q^{-1/2-iz})X_1(2iz),\end{aligned}$$
which is equivalent to
$$\begin{aligned}&
C_1(1/2 + iz)(1 + q^{1/2+iz})X_1(-2iz)\\
=&C_1(1/2 - iz)q^{4giz}(1 + q^{-1/2-iz})X_1(2iz).\end{aligned}\leqno(10)$$
Substituting the definition of $X_1$ in (10), we have
$$\begin{aligned}&
C_1(1/2 - iz)(1 + q^{-1/2-iz})q^{2giz}
\prod_{j=1}^g
(1 -\alpha_j q^{-2iz})(1 -\beta_j q^{-2iz})\\
=&C_1(1/2 + iz)(1 + q^{-1/2+iz})q^{-2giz}\prod_{
j=1}^g
(1 -\alpha_j q^{2iz})(1 -\beta_j q^{2iz}).\end{aligned}\leqno(11)$$

It suffices to derive a contradiction assuming that $z$ is not real. Put
$z = x+iy,\ x, y\in \mathbb R$. If $s = 1/2+iz$ is a zero of $\widehat\zeta_2(s)$, then $1-s = 1/2-iz$ is
also a zero. Therefore we may assume that $y > 0$. We compare the absolute
value of the both sides of (11). Put $w = q^{2iz}$. Then $|w| < 1$. By Lemma 1,
we have
$$|\prod_{
j=1}^g
(w -\alpha_j)(w -\beta_j )|\, >\,| 
\prod_{
j=1}^g
(1 -\alpha_j w)(1 -\beta_j w)|,$$ 
which implies
$$\begin{aligned}&| q^{2giz}\prod_{
j=1}^g
(1-\alpha_j  q^{-2iz})(1-\beta_j  q^{-2iz} )|\,\\
 >&\,| q^{-2giz} 
\prod_{
j=1}^g
(1 -\alpha_j  q^{2iz})(1 -\beta_j  q^{2iz})|.\end{aligned}\leqno(12)$$
Similarly we have
$$\begin{aligned}&|C_1(1/2 - iz)(1 + q^{-1/2+iz})|\\
=&|q^{a(1/2-iz)}(1 + q^{-1/2+iz})(1 + q^{-1/2-iz})q^{h(iz-1/2)}
\prod_{
j=1}^h
(1 -\gamma_j q^{-iz})(1 -\delta_jq^{-iz})|\\
>&|q^{a(1/2-iz)}(1 + q^{-1/2-iz})(1 + q^{-1/2+iz})q^{h(-iz-1/2)}
\prod_{
j=1}^h
(1 -\gamma_j q^{iz})(1 -\delta_jq^{iz})|\\
=&|C_1(1/2 + iz)(1 + q^{-1/2-iz})|.\end{aligned}$$
(Here note that $(1 + q^{-1/2+iz})(1 + q^{-1/2-iz}) \not= 0$ by the remark given in the
beginning of the proof.) This is a contradiction and completes the proof.
\vskip 0.30cm
\noindent
{\bf \S3.} As a simple choice, we drop the term $q^{-hs}\prod_{
j=1}^h
(1 -\gamma_j q^{s-1/2})(1 -\delta_jq^{s-1/2})$ and take $C_1(s) = q^{as}(1 + q^{-s})$ with a nonnegative integer $a$. Then
we have
$$\widehat\zeta_2(s) =
\frac{q^{as}(1 + q^{-s})Y (2s)}{
1 - q^{1-s}}-\frac{ 
q^{a(1-s)}q^{-s}(1 + q^{s-1})Y (2s - 1)
}{1 - q^{-s}}.\leqno(13)$$
We put $t = q^{-s}$ and are going to express $\widehat\zeta_2(s)$ as a polynomial of $t$. We have
. We have
$$\begin{aligned}Y (2s) =& q^{(1-g)/2}t^{2(g-1)}\frac{
\prod_{
i=1}^g(1 -\alpha_it^2)(1-\beta_it^2)}{(1 - t^2)(1 - qt^2)},\\
Y (2s - 1)  =& q^{3(1-g)/2}t^{2(g-1)}\frac{
\prod_{
i=1}^g(1 -q\alpha_it^2)(1-q\beta_it^2)}{(1 - qt^2)(1 - q^2t^2)}.\end{aligned}$$

Therefore we have
$$\begin{aligned}\widehat\zeta_2(s)=&q^{(1-g)/2}t^{2(g-1)}t^{-a}(1+t)\frac{
\prod_{
i=1}^g(1 -\alpha_it^2)(1-\beta_it^2)}{(1-qt)(1 - t^2)(1 - qt^2)}\\
-&q^{3(1-g)/2}t^{2(g-1)}q^at^at(1+q^{-1}t^{-1})\frac{
\prod_{
i=1}^g(1 -q\alpha_it^2)(1-q\beta_it^2)}{(1-t)(1 - qt^2)(1 - q^2t^2)}\\
=&q^{(1-g)/2}t^{2(g-1)}t^{-a}\frac{
\prod_{
i=1}^g(1 -\alpha_it^2)(1-\beta_it^2)}{(1-t)(1-qt)(1 - qt^2)}\\
-&q^{3(1-g)/2}t^{2(g-1)}q^at^a\frac{
\prod_{
i=1}^g(1 -q\alpha_it^2)(1-q\beta_it^2)}{(1-t)(1 - qt)(1 - qt^2)}.\end{aligned}$$
Now assume that $g = 1$. Put $\alpha=\alpha_1$, $\beta=\beta_1$, $c=\alpha_1+\beta_1$.  Then we have
$$\begin{aligned}&t^a(1 - t)(1 - qt)(1 - qt^2)\widehat\zeta_2(s)\\
=&[(1 - ct^2 + qt^4) - q^{a-1}t^{2a}(1 - qct^2 + q^3t^4)].\end{aligned}$$

If we choose $a = 0$, then $c$ disappears, which is unnatural. Take $a = 1$. Then
we find
$$\begin{aligned}&t(1- t)(1- qt)(1- qt^2)\widehat\zeta_2(s)\\
=&[(1- ct^2 + qt^4)- t^2(1- qct^2 + q^3t^4)]\\
=&- [q^3t^6- q(c + 1)t^4 + (c + 1)t^2- 1]\\
=&- (qt^2 + 1)(q^2t^4 + (q- c- 1)t^2 + 1)\end{aligned}$$
Using $|c|\leq q+1$, we can directly verify that all roots of this polynomial have
absolute value $1/\sqrt q$, which is the Riemann hypothesis. Therefore the most
natural choice is\footnote{When $g > 1$, the other choice of $a$, say $a = g$, may become more natural.}
$$ \widehat\zeta_2(s) =
\frac{(1 + q^s)Y (2s)}{
1- q^{1-s}}-
\frac{q^{-s}(1 + q^{1-s})Y (2s- 1)}{
1- q^{-s}}.\leqno(14)$$ 

\vskip 0.30cm
\noindent
{\bf \S4.}
 We will show that a slightly simpler function
$$\widehat\zeta_2^*(s) =
\frac{Y (2s)}{
1- q^{1-s}}-\frac{
q^{-s}Y (2s- 1)}{
1- q^{-s}},\leqno(15)$$
which can be an analogue of Wengfs zeta function, does not satisfy the Riemann
hypothesis in general. $\widehat\zeta_2^*(s)$ corresponds to the choice $C_1(s) = 1$. We
assume that $\beta_i=\overline{\alpha_i}, |\alpha_i|=\sqrt q$
for $1\leq i\leq  g$. We return to (11) and put
$w = q^{iz}$. If the Riemann hypothesis is true for $\widehat\zeta_2^*(s)$, then all roots of the
equation
$$\begin{aligned}&
(1 + q^{-1/2}w^{-1})
\prod_{
j=1}^g
(w^2- \alpha_j)(w^2- \beta_j)\\
=&(1 + q^{-1/2}w)
\prod_{
j=1}^g
(1- \alpha_j w^2)(1- \beta_j w^2)\end{aligned}\leqno(11')$$
must have absolute value 1. We consider both sides of $(11')$ as functions of $w$
in the interval $-q^{1/2} < w <-1$. Let $f(w)$ (resp. $g(w)$) denote the function
on the left-hand (right-hand) side. Since $\beta_i=\overline{\alpha_i}, |\alpha_i|=\sqrt q$, both sides are
positive. By Lemma 1, we have
$$\prod_{
j=1}^g
(w^2- \alpha_j)(w^2- \beta_j)
<
\prod_{
j=1}^g
(1- \alpha_j w^2)(1- \beta_j w^2)$$
in this interval. We replace
$\prod_{
j=1}^g
(1- \alpha_jq^{-s})(1- \beta_jq^{-s})$ by
$\prod_{
j=1}^g
(1- \alpha_j q^{-s})(1- \beta_j q^{-s})^m$ for a sufficiently large positive integer $m$.\footnote{I do not know whether such process appears in geometry of algebraic curves. Probably it does.} Then we see that
$f(w_0) < g(w_0)$ for some point $w_0 \in (-q^{1/2},-1)$. On the other hand, we
have $f(-q^{1/2}) > g(-q^{1/2}) = 0$. Therefore $f(w_1) = g(w_1)$ holds for some
$w_1 \in (-q^{1/2},w_0)$. This shows that $\widehat\zeta_2(s)$ does not satisfy the Riemann hypothesis.
\vskip 1.0cm
\centerline{\bf REFERENCES}
\vskip 0.30cm
\noindent
[K] H. Ki, All but finitely many non-trivial zeros of the approximations of the Epstein zeta function are simple and on the critical line, Proc. London  Math. Soc (3) 90 (2005), 321-344 
\vskip 0.20cm
\noindent
[LS] J. Lagarias \& M. Suzuki, The Riemann Hypotheis for certain integrals of Eisenstein series, J. Number Theory 118 (2006), 98--122
\vskip 0.20cm
\noindent
[W] L. Weng, Rank two non-abelian zeta and its zeros, preprint 2004, arXiv:math\\
/0412009, a short version is published at J. Ramanujan Math. Soc., 21(2006), 205-266.

\vskip 1.0cm
Hiroyuki Yoshida

Department of Mathematics

Faculty  of Science

Kyoto University

Kyoto 606-8502, Japan

Email: yoshida@math.kyoto-u.ac.jp
\end{document}